\documentclass{article}
\usepackage{amsmath,amsthm,amssymb}

\newcommand{\lib}{Libu{\v s}e}
\newcommand{\prem}{P{\v r}emysl}
\newcommand{\old}{Old\v rich}
\newcommand{\boz}{Bo{\v z}ena}
\newcommand{\RR}{\mathit{MRR}}

\newcommand{\fb}{\mathfrak{b}}

\newcommand{\fd}{\mathfrak{d}}

\newcommand{\last}{\mathit{last}}

\newcommand{\gw}{\omega}

\newcommand{\gS}{\Sigma}
\newcommand{\gs}{\sigma}

\newcommand{\R}{\mathbb{R}}

\newcommand{\cov}{\mathtt{cov}}

\newcommand{\cof}{\mathtt{cof}}
\newcommand{\non}{\mathtt{non}}
\newcommand{\meager}{\mathtt{meager}}
\newcommand{\lnull}{\mathtt{null}}

\newcommand{\dotrgen}{{\dot r}_{\mathit{gen}}}

\newcommand{\vecrgen}{{\vec r}_{\mathit{gen}}}

\newcommand{\cmin}{c_{\mathrm{min}}}

\newcommand{\rng}{\mathrm{rng}}

\newtheorem{theorem}{Theorem}[section]
\newtheorem{lemma}[theorem]{Lemma}
\newtheorem{claim}[theorem]{Claim}
\newtheorem{corollary}[theorem]{Corollary}
\newtheorem{fact}[theorem]{Fact}

\theoremstyle{definition}
\newtheorem{definition}[theorem]{Definition}
\newtheorem{example}[theorem]{Example}
\newtheorem{question}[theorem]{Question}

\title{Proper forcing and rectangular Ramsey theorems\footnote{2000 AMS subject classification 03E60, 03E17, 03E02}}
\author{Jind{\v r}ich Zapletal\thanks{Partially supported by grant GA {\v C}R 
201-03-0933, NSF grants DMS 0071437 and DMS 0300201, and a visiting appointment
at CRM, Universitat Aut{\` o}noma de Barcelona.}\\University of Florida}

\begin{document}

\maketitle
\begin{abstract}
I prove forcing preservation theorems for products of definable partial orders preserving cofinality
of the meager or null ideal. Rectangular Ramsey theorems for the related ideals follow from the proofs.
\end{abstract}  

\section{Introduction}
The following fact is the archetype of a rectangular Ramsey theorem.

\begin{fact}
\label{firstfact}
\cite{kechris:classical} 19.6.
For every Borel coloring $\R^2=\bigcup_n A_n$ of the plane there is a monochromatic rectangle with perfect sides.
\end{fact}

Many variations on it have been proved, for higher dimensions as well as with other demands on the size of the sides
of the monochromatic rectangle. In this paper I will present a theorem which yields many generalizations of
the above fact, connects the rectangular Ramsey theory with the theory of proper forcing and determinacy, and
isolates possible new open questions.

In order to facilitate the statement of the theorem and the surrounding discussion, it is convenient to introduce
the following notation. 

\begin{definition}
Suppose that $\vec I$ is a finite or countable sequence of $\gs$-ideals
on the reals indexed by elements of some set $X$. 
A \emph{block} is a subset of $\R^X$ of the form $\prod\vec B=\{\vec r\in\R^X:\forall x\in X\ \vec r(x)\in B_x\}$, 
where $\vec B=\langle B_x:x\in X\rangle$ is a sequence whose coordinates are Borel
sets of reals positive with respect to the corresponding ideal in the $\vec I$-sequence.
$\RR(\vec I)$ (for \emph{mutually rectangularly Ramsey})
is the statement that the collection of Borel subsets of $\R^X$ containing no block as a subset, is a $\gs$-ideal. 
If $\RR(\vec I)$ holds then $\prod\vec I$ denotes this $\gs$-ideal. If $\vec I$ is a sequence of arbitrary length,
$\RR(\vec I)$ is the statement that $\RR(\vec J)$ holds for every countable subsequence $\vec J$ of $\vec I$.
\end{definition}

\noindent Fact~\ref{firstfact} is just the statement $\RR(J, J)$, where $J$ is the $\gs$-ideal of countable sets on the reals.
Statements of the form $\RR(\vec I)$ are usually proved using sophisticated methods in Ramsey theory such as
Millikan's theorem. One point in this paper is that they can be frequently proved just from the topological forcing
properties of the factor posets $P_{I}$ for ideals $I$ in the sequence $\vec I$, where 

\begin{definition}
If $I$ is a $\gs$-ideal on the reals,
the symbol $P_I$ denotes the factor poset of Borel $I$-positive sets of reals ordered by inclusion. 
If $\vec I$ is a sequence of ideals then $P_{\vec I}$ is the side-by-side countable support product of the
factor posets of $\gs$-ideals in the sequence. Note that if the sequence $\vec I$ is countable and $\RR(\vec I)$
holds, then $P_{\vec I}$ naturally densely embeds into the poset $P_J$ where $J=\prod\vec I$.
\end{definition}

\noindent The key forcing property in this context turns out to be

\begin{definition}
$\Phi(P)$ is the statement that $P$ is proper, and every Borel meager set in the extension is covered
by one coded in the ground model.
\end{definition}

\noindent This is a rather traditional property known to be preserved under countable support iterations--\cite
{bartoszynski:set}, Theorem 6.3.22, or \cite{z:book}, Theorem 5.4.10.
I will show that in the definable context it is preserved even under countable support side-by-side products:

\begin{theorem}
\label{firsttheorem}
(LC) Suppose that $\vec I$ is a sequence of definable $\gs$-ideals, and suppose that $\Phi(P_{I})$
holds for every ideal $I$ in it. Then $\Phi(P_{\vec I})$ and $\RR(\vec I)$ hold.
\end{theorem}

The rectangular Ramsey property falls out
from the argument as a consequence of the computation of the ideal associated with the product. 
A variation of the proof gives another preservation theorem. Consider the following property of partial orders:

\begin{definition}
$\Psi(P)$ is the statement that $P$ is proper and every Borel Lebesgue null set in the extension
is covered by one coded in the ground model.
\end{definition}

\noindent This is another forcing property preserved by countable support iterations--
\cite{bartoszynski:set}, 6.3.F, or \cite{z:book}, 5.4.10. It implies $\Phi$--
\cite{bartoszynski:set}, Theorem 2.3.1, and it turns out
that it is very closely related to $\Phi$.

\begin{theorem}
\label{secondtheorem}
(LC) Suppose that $\vec I$ is a sequence of definable $\gs$-ideals on the reals.
If $\Psi(P_I)$ holds for every ideal $I$ in it. Then $\RR(\vec I)$ and $\Psi(P_{\vec I})$ hold.
\end{theorem}

A similar argument gives an asymmetric preservation result. Consider the following property of partial orders.

\begin{definition}
$\Theta(P)$ is the statement that $P$ is proper and the set of ground model reals is not meager in
the extension.
\end{definition}

\noindent Again, this forcing property is preserved under definable iterations. The preservation theorems
in the above style fail for side-by-side products of definable posets satisfying $\Theta$. However,
the following is still true.

\begin{theorem}
\label{thirdtheorem}
(LC) Suppose that $I, J$ are definable $\gs$-ideals such that $\Theta(P_I)$ and $\Phi(P_J)$ hold.
Then $\RR(I,J)$ and $\Theta(P_I\times P_J)$ hold.
\end{theorem}

The theorems in this paper work only in definable context and use the appropriate large cardinal assumptions
(denoted by LC) to guarantee the determinacy of certain infinite games with real entries. One way to make this
precise is to demand that for the relevant ideals the collection of (codes for) analytic sets which belong to the ideal
is projective, assume that there are $\gw+\gw$ Woodin cardinals, and use the results of forthcoming \cite{neeman:book}.
The notation used in the paper sticks to the set theoretic standard of \cite{jech:set}. For a finite
binary sequence $s$ the term $[s]$ denotes the open set of all infinite binary sequences containing
$s$ as an initial segment.  

\section{Several examples}

While there are not many definable partial orders $P$ with $\Phi(P)$ that occured in practice, they do not seem
to admit a simple classification. This is partially due to the fact that the property is preserved by countable
support iterations and products. The following examples are more or less well-known \cite{z:book}.

\begin{example}
Let $I$ be the ideal of countable subsets of $\R$. The Sacks forcing can be densely embedded in the factor $P_I$,
and Fact~\ref{firstfact} just translates into the statement $\RR(I,I)$.
\end{example}

\begin{example}
Let $E_0$ be the modulo finite equality on $2^\gw$ and let $I$ be the $\gs$-ideal generated by Borel sets
which visit each $E_0$-equivalence class in at most one point.The rectangular Ramsey property 
$\RR(I, I, I\dots)$ is proved for example in \cite{z:book} 2.3.16,
which also gives a computation a simple basis for the ideal $\prod(I, I, I,\dots)$ in terms of the ideal $I$.
\end{example}

\begin{example}
Let $G$ be the graph on $2^\gw$ connecting two binary sequences just in case they differ in exactly one point.
Let $I$ be the $\gs$-ideal generated by the Borel $G$-independent subsets of $2^\gw$. The Silver forcing densely
embeds into $P_I$. The situation is parallel to that of $E_0$-forcing.
\end{example}

\begin{example}
Let $\cmin$ be the clopen partition of pairs of infinite binary sequences in two classes defined by
$\cmin(x, y)=x\Delta y\mod 2$ where $x\Delta y$ is the smallest number where the sequences $x, y$ differ. Let
$I$ be the $\gs$-ideal generated by $\cmin$ homogeneous sets. $P_I$ is a forcing similar to Sacks forcing, and it
satisfies $\Phi(P_I)$ as well.
\end{example}

All of these forcings actually satisfy the stronger property $\Psi$. It is not an entirely trivial matter
to find a poset satisfying $\Phi$ but not $\Psi$; this is done for example in Section 7.3.C of
\cite{bartoszynski:set}.

The non-examples are perhaps more interesting than the examples. Exactly which instances
of rectangular Ramsey properties does Theorem~\ref{firsttheorem} not cover
and why? Recall the relevant part of Cicho{\' n}'s diagram:

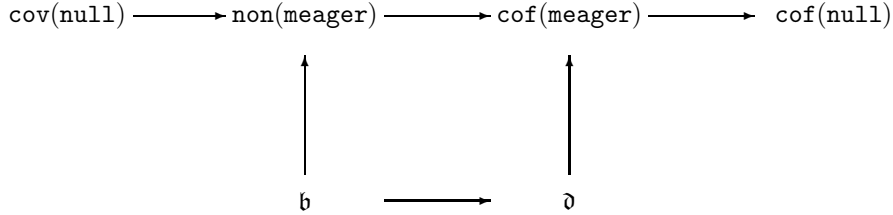
\begin{picture}(330,120)
 \put(35,100){\vector(1,0){35}}
 \put(130,100){\vector(1,0){40}}
 \put(230,100){\vector(1,0){40}}
 \put(200,40){\vector(0,1){45}}
 \put(100,40){\vector(0,1){45}}
 \put(130,30){\vector(1,0){40}}
 \put(100,30){\makebox(0,0){$\fb$}}
 \put(10,100){\makebox(0,0){$\cov(\lnull)$}}
 \put(100,100){\makebox(0,0){$\non(\meager)$}}
 \put(200,100){\makebox(0,0){$\cof(\meager)$}}
 \put(300,100){\makebox(0,0){$\cof(\lnull)$}}
 \put(200,30){\makebox(0,0){$\mathfrak{d}$}}
\end{picture}

\noindent It is well-known that $\cof(\meager)$ is the maximum of $\fd$ and $\non(\meager)$--\cite{bartoszynski:set},
Theorem 2.2.11. Thus forcings which fail $\Phi$ either must fail to be bounding or must make the
ground model reals meager.  What happens if the ideals in question violate the bounding condition?

\begin{example}
Let $I$ be the $\gs$-ideal generated by compact subsets of $\gw^\gw$. The Miller forcing densely embeds into
the factor $P_I$. Spinas \cite{spinas:ramsey} proved that $P_I\times P_I$ is proper and in fact $\RR(I, I)$ holds. However,
$\RR(I, I, I)$ fails rather badly by a result of \cite{velickovic:complexity}: there is a Borel
map $f:(\gw^\gw)^3\to\gw^\gw$ such that the image of every block with superperfect sides contains a nonempty open set.
Then for every dense codense set $X\subset\gw^\gw$ the coloring $c:(\gw^\gw)^3\to 2$ given by $
c(x,y,z)=0$ iff $f(x,y,z)\in X$
witnesses the failure of $\RR(I,I,I)$. 

It is now in fact easy to observe
that if $J_i:i\in 3$ are $\gs$-ideals such that $P_{J_i}$ are proper not bounding notions of forcing, then
$\lnot\RR(J_0, J_1, J_2)$.
Fix a Borel coloring $c:(\gw^\gw)^3\to 2$ with no monochromatic block with superperfect sides. The assumption
on the ideals $J_i$ implies that there are $J_i$-positive Borel sets $B_i$ and Borel functions $f_i:B_i\to\gw^\gw$
such that for every $J_i$-positive Borel subset $C_i\subset B_i$ the set $f_i''C_i$ is not $\gs$-compact.
Let $d:\prod_{i\in 3} B_i\to 2$ be defined by $d(x, y, z)=c(f_0(x), f_1(y), f_2(z))$. It is clear from the construction that
there is no monochromatic block with $J_i$-positive Borel sides.
\end{example}

\noindent Now what if the poset $P_I$ violates the nonmeagerness condition?

\begin{example}
Consider the ideal $I$ of Lebesgue null sets of the reals with its attendant
Solovay forcing $P_I$. Arnold Miller pointed out to me the most elegant way to see that $\RR(I, I)$ fails. 
It is a well-known fact 
(Theorem 3.2.10 of \cite{bartoszynski:set}) that for every two Borel $I$-positive sets
$A, B$ the set $A+B$ contains an interval. Then argue just like in the previous example.
\end{example}

\noindent Looking back at Cicho{\' n}'s diagram, it is clear that we are squeezed into a quite small space in the search for 
ideals $I$ for which $\RR(I,I,I)$ holds but
is not implied by Theorem~\ref{firsttheorem}. 

\begin{example}
A typical bounding forcing making the ground
model reals meager without adding random reals is the fat tree forcing. A fat tree is a tree $T\subset\gw^{<\gw}$
such that for every $n\in\gw$ there is $m\in\gw$ such that every node in the tree of length at least $m$ has at least
$n$ many immediate successors. The collection $I$ of all Borel subsets of $\gw^\gw$ containing no subset of the
form $[T]$ for some fat tree $T$ is a $\gs$-ideal. I do not know whether $\RR(I,I)$, $\RR(I, I, I)$ hold.
\end{example}

\noindent The rectangular Ramsey properties concerning several different ideals are somewhat more complex.

\begin{example}
Let $I$ be the Lebesgue null ideal on $2^\gw$ and $J$ the Laver ideal on $\gw^\gw$. The Laver forcing $P_J$
is connected with the invariant $\fb$--see the above diagram and compare this case with the previous ones.
$\RR(I,J)$ fails. Consider the Borel function $f:2^\gw\times\gw^\gw\to 2^\gw$ given by $f(a,b)=a\circ b$.
It turns out that the image of every Borel rectangle with positive sides contains a
nonempty open subset of $2^\gw$. The argument is completed just as in the previous cases. A short
inspection of the proof will show that whenever $K_0,K_1$ are $\gs$-ideals such that both 
factor forcings $P_{K_0}, P_{K_1}$ are proper and one of them adds a splitting real and the other a dominating real,
$\RR(K_0,K_1)$ must fail. Note that the corresponding covering numbers $\fb$ and $\mathfrak{r}$ are both
provably $\leq\cof(\meager)$.
\end{example}

\begin{example}
Let $P$ be the eventually different real forcing \cite{bartoszynski:set}, 7.4.B, and let $I$ be its associated
$\gs$-ideal on $\gw^\gw$, so that $P$ is in the forcing sense equivalent to $P_I$. Let $J$ be the ideal of countable sets.
$\RR(I,J)$ fails. The easiest way to see that is to choose a perfect subset $X\subset\gw^\gw$
consisting of mutually eventually different functions, and define a function $f:\gw^\gw\times X\to 2^\gw$
by letting $f(y,x)$ to be the sequence of parities of elements of the set $\{n\in\gw:x(n)=y(n)\}$.
It turns out that the $f$-image of every rectangle with Borel $I$-positive and $J$-positive sides respectively, 
contains a nonempty open subset of $2^\gw$. Just like in the
previous arguments this means that $\RR(I,J)$ fails and in fact
$\RR(I,K)$ fails for every nonprincipal $\gs$-ideal $K$. Note that $\cov(I)\leq\non(\meager)\leq\cof(\meager)$.
\end{example} 

The forcing preservation theorems stated in the paper also bring up more general questions:

\begin{question}
Suppose that $I,J$ are definable $\gs$-ideals such that $\cof(\meager)$ is less than both $\cov(I\restriction B)$
and $\cov(J\restriction C)$ for all Borel $I$-positive sets $B$ and Borel $J$-positive sets $C$.
Does $\RR(I,J)$ hold?
\end{question}

\begin{question}
Are there (necessarily undefinable) $\gs$-ideals $I,J$ such that both $\Phi(P_I)$ and $\Phi(P_J)$ hold,
but $\RR(I,J)$ fails?
\end{question}

\section{The fusion games}

Fix a definable $\gs$-ideal $I$, for simplicity and without loss of generality assume that its underlying
space is $2^\gw$. Consider the following infinite game $G_\Phi(I)$ of length $\gw$ with real entries, between players
\old\  and \boz\ \cite{jirasek:povesti}. All moves in it are (codes for) Borel $I$-positive sets.

First, \boz\  indicates an initial set $B$. The game then has infinitely many rounds. In round $i\in\gw$, \old\ 
indicates finitely many sets $B(i, j): j\in j(i)$, and \boz\  responds to each by playing its subset $C(i, j)$.
The order of moves is $B(i, 0), C(i, 0), B(i, 1), C(i,1)\dots$ It is \old\  who calls a stop to the round $i$
at some stage $j(i)$. \old\  wins if no round dragged on for infinitely many moves, and the set $X=B\cap
\bigcap_{i\in\gw}\bigcup_{j\in j(i)}C(i, j)$ does not belong to the ideal $I$. The set $X$ will be referred to as
the \emph{result} of the play of the game $G_\Phi(I)$.

\begin{example}
\old\  has a winning strategy in the game $G_\Phi(I)$ where $I$ denotes the $\gs$-ideal of countable sets.
During the play he will stop $i$-th round after $2^i$ many moves, construct binary sequences $s(i,j)$
such that $s(i+1, 2j)$ and $s(i+1, 2j+1)$ are incompatible extensions of $s(i,j)$ and the moves
$B(i+1,2j)\subset [s(i+1, 2j)]$ and $B(i+1,2j+1)\subset [s(i+1, 2j+1)]$ will be perfect subsets of $C(i,j)$. 
In the end of the play, the resulting set will be equal to the perfect set
$\{r\in 2^\gw:\forall i\in\gw\ \exists j\in 2^i\ s(i,j)\subset r\}$.
\end{example}

\begin{example}
\boz\  has a winning strategy in the game $G_\Phi(I)$ where $I$ denotes the ideal of Lebesgue null sets. She just
plays so that the Lebesgue measure of the set $C(i, j)$ is less than $2^{-ij}$. Then clearly the result $X$
of the play is Lebesgue null, and she wins.
\end{example}

Perhaps some remarks on the nature
of the game are in order. It is clear that in a play of the game \old\  attempts to produce some sort of a fusion sequence. 
It is in general impossible for him to predict how long the different rounds will take. This would lead
to a related game which is connected to the invariant $\cof(\lnull)$ as opposed to $\cof(\meager)$. It is also
impossible for him to expose all his moves in a given round at the outset of the round. This would allow him to construct 
Borel positive sets of mutually generic reals, which cannot be done in the case of $E_0$ forcing or Silver forcing.
It is clear that for both sides a smaller move is a better move, and so the moves can be restricted to
an arbitrary dense subset of the poset $P_I$.

The following lemma records the key connection between the property $\Phi$ and the game $G_\Phi$.

\begin{lemma}
\label{keylemma}
(LC) $\Phi(P_I)$ iff \old\  has a winning strategy in the game $G_\Phi(I)$.
\end{lemma}

The right-to-left direction is easy. Fix a winning strategy $\gs$ for \old. To prove the properness of
the poset $P_I$, let $M$ be a countable elementary submodel of some large structure containing $\gs$
and let $B\in P_I\cap M$ be a condition. To find a master condition for the model $M$ below the set $B$,
let $\{ D_i:i\in\gw\}$ enumerate the open dense subsets of the poset $P_I$ in the model $M$ and simulate a play
of the game $G(I)$ against the strategy $\gs$ in which \boz\  plays the initial set $B$ and then chooses her
moves $C(i,j)$ to come from the sets $M\cap D_i$. It is easy to argue inductively that this is possible and all
the moves of the play will be in the model $M$. The result $X$ of the play will be an $I$-positive Borel set
below the condition $B$, and clearly the required master condition for the model $M$.

To see that every Borel meager set in the extension is covered by one coded in the ground model, let
$B\in P_I$ be a condition and $\{\dot O_n:n\in\gw\}$ a name for a sequence of open dense subsets
of $2^{<\gw}$. Fix a bijection $f:\gw\to 2^{<\gw}\times\gw$ and simulate a run of the game $G_\Phi$
in which \old\ follows his strategy $\gs$, \boz\ plays the initial set $B$ and in round $i\in\gw$
plays sets $C(i,j)$ in such a way that there is a descending chain $\langle t(i,j):j\leq j(i)\rangle$
of binary sequences such that writing $f(i)=\langle s,n\rangle$, the chain begins with $t(i,0)=s$
and $C(i,j)\Vdash \check t(i, j+1)\in\dot O_n$. This is easily possible. Clearly, the resulting set $X\subset B$
of the play will force that $t(i,j(i))\in O_n$. Thus the set $N_n=\{t\in 2^{<\gw}:X\Vdash\check t\in\dot O_n\}\subset 2^{<\gw}$
is open dense for every number $n\in\gw$, and the condition $X$ forces
the meager set $\{r\in 2^\gw:\exists n\ \forall m\ r\restriction m
\notin O_n\}$ in the extension to be a subset of the meager set $\{r\in 2^\gw:\exists n\ \forall m\ r\restriction m
\notin N_n\}$ coded in the ground model.

The opposite direction of Lemma~\ref{keylemma} is harder. The key tool is the
reduction of the game $G_\Phi(I)$ to an integer game $H_\Phi(I,C)$, where $C$ is a compact $I$-positive set. This is
a game of length $\gw$ between players \prem\  and \lib\ \cite{jirasek:povesti}. 
The game has infinitely many rounds. In each round $i$
\prem\  produces a finite collection $k(i,j):j\in j(i)$ of natural numbers, and \lib\ answers each with a
sequence $s(i, j)\in 2^{k(i,j)}$. The order of moves is $k(i,0),s(i,0),k(i,1),s(i,1),\dots$ It is \lib\ who decides
to call a stop to the round $i$ at some stage $j(i)$. \lib\ wins if 
no round dragged on for infinitely many moves and the set $X=\{r\in C:\forall i\in\gw\exists j\in j(i)\ 
s(i,j)\subset r\}$ is $I$-positive. The set $X$ will be referred to as the \emph{result} of the play.

\begin{claim}
\label{integerclaim}
(LC) $\Phi(I)$ implies that \lib\ has a winning strategy in the game $H_\Phi(I, C)$ for every compact $I$-positive set $C$.
\end{claim}

\begin{proof}
The game is determined by our large cardinal assumptions. So it is enough to derive a contradiction
from the assumption that \prem\  has a winning strategy $\gs$ for some game $H_\Phi(I, C)$.

First observe that in this case \prem\  has a \emph{positional} winning strategy $\tau$ in the form of an increasing sequence
$\langle k_l:l\in\gw\rangle$ of natural numbers such that if he plays them successively he wins no matter what
\lib's moves are. This follows from the fact that a larger number is a better move, as far as \prem\  is concerned.
To define the numbers $k_l$ note that at each move of the game $H(I,C)$ \lib\ has just finitely many options
at her disposal, and so for every number $l$ there is just a finite set $Z_l$ 
of numbers that the strategy $\gs$ can possibly use at the $l$-th move. Let then $k_l=\max(Z_l\cup\{k_{l-1}\})+1$.
For every play $x$ against the strategy $\tau$ with \lib's moves enumerated as $s_l$ consider the play
$y$ against the strategy $\gs$ in which \lib\ calls stops to rounds at the same places as in the play $x$ and she
plays the sequences $s_l\restriction k$ where $k\in k_l$ is the
number the strategy $\gs$ produces at the $l$-th move. Clearly, the result of the play $x$ is a subset of the
result of the play $y$, and since the strategy $\gs$ was winning, both of these resulting sets must belong
to the ideal $I$. Thus the strategy $\tau$ is winning as well.

Now fix a positional winning strategy $\tau=\langle k_l:l\in\gw\rangle$ for \prem\  as in the previous paragraph.
Let $r\in C$ be a $V$-generic real for the poset $P_I$ below the condition $C$. Let $g\in\prod_l 2^{k_l}$ be a function
in the extension defined by $g(l)=r\restriction k_l$. Since the ground model reals are not meager,
there is a ground model function $h\in\prod_l 2^{k_l}$ such that the set $u=\{l\in\gw:h(l)=g(l)\}$ is infinite.
Since the ground model reals are dominating, there is a ground model infinite set $v\subset\gw$ such that between
every two successive elements of it there is an element of the set $u$. Now consider the play of the
game $H(I,C)$ against the strategy $\tau$ in which \lib\ plays the sequences $s_l=h(l)$ and calls stops to
the rounds after each $l$-th move where $l\in v$. It follows from the definition of the resulting set $X$ of the play
and the choice of $h$ and $v$ that $r\in X$. However, this play is in the ground model, so its result $X$
is an $I$-small closed set coded in the ground model, and the generic real $r$ is forced to fall out of all such sets.
Contradiction!
\end{proof}

For the reduction of the game $G_\Phi$ to the game $H_\Phi$ I will need a small general observation which does
not concern these games.

\begin{claim}
If the poset $P_I$ is proper and bounding, then compact sets are dense in it. Moreover, for every countable
elementary submodel $M$ of a large structure and a condition $B\in M\cap P_I$, there is a compact $I$-positive set
$C\subset B$ such that for every dense set $D\subset P_I$ in the model $M$ there is a number $k\in\gw$
such that for every sequence $s\in 2^k$ there is an element $E\in D\cap M$ such that $C\cap [s]\subset A$.
\end{claim}

\begin{proof}
The first sentence is proved in \cite {z:book}. For the rest of the claim, fix the model $M$ and a condition $B\in M\cap
P_I$. Enumerate the infinite maximal antichains in $M$ consisting of compact sets by $\{A_n:n\in\gw\}$, and for each of
them enumerate the set $A_n\cap M$ by $\{E(n, m):m\in\gw\}$. Finally, choose a master condition
$B'\subset B$ for the model $M$. 

Note that $B'\Vdash\forall n\ \exists!m\ E(n,m)\in\dot G$, and let $\dot f\in\gw^\gw$
be the name for a function assigning to each number $n$ the unique $m$ such that $E(n,m)\in\dot G$. The forcing
$P_I$ is bounding, and so there is a ground model function $g\in\gw^\gw$ and a condition $B''\subset B'$
such that $B''\Vdash\dot f<\check g$ pointwise. Let $B'''=B'\cap\bigcap_{n\in\gw}\bigcup_{m\in g(n)}E(n,m)$.
The set $B'''$ is Borel, and since the condition $B''$ forces the generic real into it, it has to be $I$-positive.
Let $C\subset B'''$ be any compact $I$-positive subset.

Now clearly for every number $n\in\gw$ the finitely many compact sets $\{E(n,m):m\in g(n)\}$ in the antichain $A_n$
cover the compact set $C$, and by a compactness argument there must be a number $k$ such that for every sequence $s\in 2^k$
there is some $m\in g(n)$ such that $C\cap [s]\subset E(n,m)$. The claim immediately follows.
\end{proof}

Now assume that $\Phi(P_I)$ holds. Note that the large cardinal assumptions imply that the game $G_\Phi(I)$
is determined, and it is enough to obtain a contradiction from the assumption that \boz\  has a winning strategy
$\gs$ in it. Let $M$ be a countable elementary submodel of a large structure containing the strategy $\gs$
and let $B\in M\cap P_I$ be the initial move dictated by the strategy. Let $C\subset B$ be the compact
set from the previous claim, and let $\tau$ be \lib's winning strategy in the game $H(I, C)$ from Claim~\ref{integerclaim}.
The contradiction will be reached by pitting the strategies $\tau$ and $\gs$ against each other in a way.

Find plays $x$ and $y$ of the games $G_\Phi(I)$ and $H_\Phi(I, C)$ observing the strategies $\gs$ and $\tau$ respectively
so that

\begin{itemize}
\item all moves of the play $x$ are in the model $M$
\item \old\  calls stops to rounds in the play $x$ exactly when \lib\ calls stops in the play $y$
\item at move $l$, writing $C_l$ for the $l$-th \boz's move in the play $x$
and $s_l$ for the $l$-th \lib's move in the play $y$ it is the case
that $C\cap O_{s_l}\subset C_l$.
\end{itemize}

After this is done, it is clear from the second and third items that the result of the play $y$ is a subset
of the result of the play $x$, and since \lib's strategy $\gs$ was winning, it must be the case that
both of the resulting sets must be $I$-positive. This means that \boz\  lost the play $x$ observing her strategy
$\gs$, and this contradiction will finish the proof.

The plays $x$ and $y$ are built simultaneously by induction. Suppose that $l\in\gw$ is a number and
the partial plays $x\restriction l$
and $y\restriction l$ have been built. The set $Z=\{E\in P_I:$ there is $D\in P_I$ such that
$x\restriction l^\smallfrown D^\smallfrown E$ is a play of the game $G_\Phi(I)$ observing the strategy $\gs\}$
is dense in $P_I$ and belongs to the model $M$. By the choice of the set $C$ there is a number $k_l\in\gw$
such that for every sequence $s\in 2^{k_l}$ there is a set $E\in Z\cap M$ such that $C\cap [s]\subset E$.
Put $k_l$ to be the next \prem's move in the play $y$. The strategy $\tau$ responds with some sequence
$s_l$. By the choice of the number $k_l$ there are sets $B_l$ and $C_l$ in the model $M$
such that the play $x\restriction(l+1)=x\restriction l^\smallfrown B_l^\smallfrown C_l$ 
respects the strategy $\gs$ and $C\cap [s_l]\subset C_l$. This completes the induction step
and the proof of Lemma~\ref{keylemma}.

\section{The products}
\label{productsection}

Theorem~\ref{firsttheorem} now follows by a rather standard, if notationally awkward, fusion argument. 
For simplicity I will treat only the case of a countable sequence $\vec I
=\langle I_n:n\in\gw\rangle$ of definable
$\gs$-ideals. Suppose $\Phi(P_{I_n})$ holds for every number $n\in\gw$. Define a game $G$ between \old\  and \boz\ 
in the same way as the game $G_\Phi(I)$ was defined, only the
moves now will be blocks ordered by inclusion and
\old\  wins if the result of the play contains a block. To simplify the notation, the initial block played
by \boz\ will be denoted by $p_{ini}$, and then the blocks played by \old\ and \boz\ will be indexed
as $\{p_k:k\in\gw\}$ and $\{q_k:k\in\gw\}$ respectively, in the increasing order--so $q_k\subset p_k$.
The $n$-th coordinate of the blocks $p_k, q_k$ will be denoted by $p_k(n), q_k(n)$ respectively.

\begin{claim}
\old\  has a winning strategy in the game $G$.
\end{claim}

\begin{proof}
Choose winning strategies $\{\gs_n:n\in\gw\}$ for \old\  in the games $\{G_\Phi(I_n):n\in\gw\}$. The winning strategy $\gs$
in the game $G$ will be their fusion of sorts. It is fully determined by the following demands:

\begin{itemize}
\item on the side, the strategy $\gs$ will inductively obtain infinite sets $\gw=u_0\supset u_1\supset\dots$; the set
$u_n$ will be enumerated in the increasing order as $\{l(m,n):m\in\gw\}$.
\item for all numbers $n$ and $m$, $p_{l(m,n)}(n)\supset q_{l(m,n)}(n)=
p_{l(m,n)+1}(n)\supset q_{l(m,n)+1}(n)=p_{l(m,n+2)}(n)\supset\dots
\supset q_{l(m+1,n)-1}(n)$. Moreover, $p_{ini}(n)=p_0(n)\supset q_0(n)=p_1(n)\supset\dots\supset q_{l(0,n)-1}(n)$.
This tells the strategy $\gs$ how to determine the $n$-th coordinate of the condition
$p_k$ at the $k$-th move if $k\notin u_n$.
\item for all numbers $n$, the play $x_n$ of the game $G_\Phi(I_n)$ given by the initial move $q_{l(0,n)-1}$ and then the
exchange $p_{l(0,n)}(n), q_{l(1,n)-1}(n), p_{l(1,n)}(n), q_{l(2,n)-1}(n)\dots$, respects the strategy $\gs_n$.
This tells the strategy $\gs$ how to determine the $n$-th coordinate of the condition $p_k$ at move $k$ if $k\in u_n$.
Note that the previous item implies that \boz's answers in this play are legal.
\item $u_{n+1}=\{k\in u_n:$ the strategy $\gs_n$ calls a stop to a round in the play $x_n$ just before the move $p_k(n)\}$.
This tells the strategy $\gs$ how to build the set $u_{n+1}$.
\item the stops to the rounds will be called before each stage indexed by $\min(u_n)$ for some $n\in\gw$.
\end{itemize}

In other words, at $0$-th coordinate \old\ just simulates a play $x_0$ of the game $G_\Phi(I_0)$ respecting the
strategy $\gs_0$. He puts $u_1$ to be the set of indexes of moves before which stops are called in $x_0$,
and at $1$-th coordinate at moves indexed by numbers in the set $u_1$
he will simulate a play $x_1$ of the game $G_\Phi(I_1)$ respecting the strategy
$\gs_1$. At other moves he just plays dead--repeats the last \boz's move. Similarly for the other coordinates.

It must be proved that $\gs$ is a winning strategy. For every number $n\in\gw$, look at the play $x_n$ of the game 
$G_\Phi(I_n)$ 
from the third item above and let $X_n\subset 2^\gw$ be its resulting set. Since the strategy $\gs_n$ is winning for \old,
this set is Borel and $I$-positive, and it will be enough to show that the result of the whole play of the game
$G$ contains the block $\prod_{n\in\gw} X_n$.

So let $\vec r\in\prod_n X_n$ be an $\gw$-sequence, and let $m\in\gw$ be an arbitrary number. I must produce
a number $j$ with $\min(u_n)<j\leq\min(u_n)$ such that $\vec r\in q_{j-1}$. First note that for all numbers
$n\geq m+1$ it is the case that $\vec r(n)\in q_{j-1}(n)$ for all numbers $j$ in this range, since the set $X_n$
is a subset of the initial move of the play $x_n$ which in turn is a subset of the set $q_{j-1}(n)$ by the second
item above. To handle the numbers $m\leq n$, by downward induction on $m\leq n$ construct a nondecreasing sequence
$\{j(m):m\leq n\}$ of natural numbers such that

\begin{itemize}
\item $j(n)\in u_n$ so that $\min(u_n)<j(n)\leq \min(u_{n+1})$.
Moreover, for every number $m<n$, $j(m)\in u_m$ and $\max(u_{m+1}\cap j(m+1))<j(m)\leq j(m+1)$
\item $\vec r(m)\in q_{j(m)-1}(m)$.
\end{itemize}

This is easily possible to arrange--to find the number $j(m)$ note that the set $X_m$ is the result of the
play $x_m$ and $\vec r(m)\in X_m$. It is not difficult to verify from the definition of the strategy $\gs$
that the  number $j=j(0)$ is as required.
\end{proof}

\begin{corollary}
The collection of Borel subsets of $(2^\gw)^\gw$ containing no block is a $\gs$-ideal.
\end{corollary}

\begin{proof}
Let $p_0=\prod_{n\in\gw}B_n$ be a block, decomposed into a countable union $p_0=\bigcup_mC_m$ of Borel sets.
It is enough to show that one of the sets $C_m$ contains a block.

Write $\vecrgen$ for the $P_{\vec I}$-generic sequence of reals. Note that $p_0\Vdash\vecrgen\in\dot p_0$,
since for every number $n\in\gw$ the $n$-th coordinate $\vecrgen(n)$ is forced to be $P_{I_n}$-generic
below the set $B_n$ and therefore to belong to the set $B_n$ in the extension. By an absoluteness argument,
there is a block $p_1\subset p_0$ which forces $\vecrgen\in\dot C_m$ for some definite number $m\in\gw$.
I claim that the set $C_m$ contains a block.
 
To see this, let $M$ be a countable elementary submodel of a large enough structure containing the
condition $p_1$, the set $C_m$, as well as a winning strategy $\gs$ for \old\  in the game $G$. Enumerate
all open dense subsets of the poset $P_{\vec I}$ in the model $M$ as $\{D_i:i\in\gw\}$ and simulate
a run $x$ of the game $G$ against the strategy $\gs$ with the initial move $p_1$ and such that
all its moves are in the model $M$ and during the $i$-th round \boz\  plays only sets from
the open dense set $D_i$. The result of the game contains some block
$p_2\subset p_1$. I claim that $p_2\subset C_m$; this will complete the proof.

Let $\vec r$ be a sequence from the block $p_2$. It is easy to check that the collection
of all blocks in the  model $M$
containing the sequence $\vec r$ is a filter on the poset $P_{\vec I}\cap M$. By the simulation above,
this filter is $M$-generic and $\vec r$ is its associated generic real. By the forcing theorem
applied in the model $M$, $M[\vec r]\models\vec r\in C_m$, and by an absoluteness argument $\vec r\in C_m$.
So $p_2\subset C_m$ as desired.
\end{proof}

\begin{corollary}
$\Phi(P_{\vec I})$ holds.
\end{corollary}

\begin{proof}
Writing
$J=\prod\vec I$ it is now 
clear that the poset $P_{\vec I}$ naturally densely embeds into $P_J$,
the game $G$ is just $G_\Phi(J)$ under another name, and \old\  has a winning strategy in it.
A reference to Lemma~\ref{keylemma} concludes the argument.
\end{proof}

This completes the proof of Theorem~\ref{firsttheorem}. I will state two corollaries of independent interest.
To facilitate the notation in the statement and proof, for a given number
$n$ let $\vec I\ominus n$ be the sequence $\vec I$ with the $n$-th entry removed. Clearly
$P_{\vec I}=P_{\vec I\ominus n}\times P_{I_n}$. Let $\vecrgen$ be the $P_{\vec I}$-generic sequence,
and let $\vecrgen\ominus n$ be just  $\vecrgen$ with its $n$-th
entry removed, understood now as the $P_{\vec I\ominus n}$-generic sequence.

\begin{corollary}
For every number $n\in\gw$, $P_{\vec I}\Vdash\vecrgen(\check n)$ belongs to no Borel $\dot I_n$-small
set coded in the model $V[\vecrgen\ominus\check n]$. 
\end{corollary}

Of course, the real $\vecrgen(n)$ belongs to no Borel $I_n$-small set coded in $V$, since it is $V$-generic
for the poset $P_{I_n}$. The point of the corollary is that the real falls out even from all $I_n$-small Borel sets
coded in the larger model $V[\vecrgen\ominus n]$.

\begin{proof}
Suppose $p$ is a $(\vec I\ominus n)$-block, $\dot U$ be a $P_{\vec I\ominus n}$-name
for an $I_n$-small set, and let $B\in P_{I_n}$ be a Borel set. It will be enough to find a $(\vec I\ominus n)$-block
$q\subset p$ and a Borel $I_n$-positive set $C\subset B$ such that $\langle q, C\rangle\Vdash\dot r\notin\dot U$
in the product $P_{\vec I\ominus n}\times P_{I_n}$, where $\dot r$ is the name for the $P_{I_n}$-generic real.

Thinning out the block $p$ if necessary we may assume that there is a Borel set
$X\subset p\times B$ such that all vertical sections of the set $X$ are $I_n$-small
and $p\Vdash\dot U\cap\dot B$ is the vertical section of the set $\dot X$ corresponding to
the sequence $(\vecrgen\ominus n)\in\dot p$. Theorem~\ref{firsttheorem} now implies 
$\RR(\prod\vec I\ominus n, I_n)$, and so either
the set $X$ or its complement in $p\times B$ must contain a rectangle with positive sides. Well, it cannot be the
set $X$ since its vertical sections are $I_n$-small, so there must be a rectangle $q\times C\subset (p\times B)\setminus X$.
It is immediate that $q,C$ work as required. 
\end{proof}

\begin{corollary}
If every ideal on the sequence $\vec I$ is $\bf\Pi^1_1$ on $\bf\gS^1_1$ then so is the ideal $\prod\vec I$.
\end{corollary}

\begin{proof}
It follows from Lemma C.0.9 of \cite{z:book} that if $J$ is a $\gs$-ideal such that the factor poset is
proper and bounding, then $J$ is $\bf\Pi^1_1$ on $\bf\gS^1_1$ iff the set of $J$-positive compact sets
is analytic iff there is an analytic dense collection of $J$-positive compact sets. Now write $J=\prod\vec I$.
Theorem~\ref{firsttheorem} implies that the factor forcing is proper and bounding.
If every ideal on the sequence $\vec I$ is $\bf\Pi^1_1$ on $\bf\gS^1_1$, then clearly the collection
of all blocks of the form $\prod\vec C$, where $\vec C$ is a sequence of compact sets positive with respect
to the corresponding ideal on the sequence $\vec I$, is analytic and dense in $P_J$. The corollary follows.
\end{proof}

\section{Cofinality of the null ideal}

The proof of Theorem~\ref{secondtheorem} is almost identical. It uses the following key combinatorial fact:

\begin{fact}
\label{nullfact}
\cite{bartoszynski:set}, Section 2.3.
$\Psi(P)$ is equivalent to properness of $P$ together with the statement ``for every ground model nondecreasing function
$h\in\gw^\gw$ diverging to infinity and for every function $f\in\gw^\gw$ in the extension, there is
a ground model function $g:\gw\to[\gw]^{<\gw}$ such that for every number $n$ the set $g(n)$ has size
at most $h(n)+1$ and contains the value $f(n)$''.
\end{fact}

The only change in the proofs is that I must devise new games $G_\Psi(I)$, $H_\Psi(I,C)$. The game $G_\Psi(I)$
is played exactly as $G_\Phi(I)$ except it is now \boz\ who decides the lengths
of the rounds, and it is her responsibility to
see to it that the lengths of the rounds are finite, never decrease and diverge to infinity.
The change in the definition of the $H$ game is the same.
Similarly as in the $\Phi$ case, the following claims are crucial:

\begin{claim}
(LC) $\Psi(P_I)$ if and only if \old\ has a winning strategy in the game $G_\Psi(I)$.
\end{claim}

\begin{claim}
(LC) $\Psi(P_I)$ implies that \lib\ has a winning strategy in the game $H_\Psi(I,C)$, for every
closed $I$-positive set $C$.
\end{claim}

\noindent The only real difference occurs in the proof of the latter claim:

\begin{proof}
The game is determined, and it will be enough to derive a contradiction from the assumption that \prem\ has
a winning strategy $\gs$. First, use a compactness argument to find a nondecreasing function $h\in\gw^\gw$
diverging to infinity such that in all plays in which \prem\ uses his strategy, the $i$-th round
will have at least $h(i)+1$ many moves.
As in Claim~\ref{integerclaim}, it is also possible to find a positional strategy $\tau$
in the form of an increasing infinite sequence $\langle k_n:n\in\gw\rangle$ of natural numbers such that \prem\ wins if he plays
an arbitrary increasing subsequence of it and lets the $i$-th round
last $h(i)+1$ steps, disregarding \lib's moves entirely.

Now let $r\in C$ be a $V$-generic real under the condition $C$, and let $e\in\prod_n 2^{k_n}$ be the function
defined by $e(n)=r\restriction k_n$. Since $\Psi(P)$ implies $\Phi(P)$, the ground model reals are not meager,
and there is a ground model function $d\in\prod_n 2^{k_n}$ such that the set $u=\{n\in\gw:e(n)=d(n)\}$
is infinite. Using Fact~\ref{nullfact} it is not difficult to find a ground model function $g:\gw\to[\gw]^{<\gw}$
such that for every number $i\in\gw$ the set $g(i)$ has size $h(i)+1$, contains some element of the set
$u$, and moreover $\max(g(i))<\min(g(i+1))$.

Consider the play of the game $H_\Psi(I,C)$ in which \prem\ plays numbers from the set
$\{k_n:n\in\bigcup\rng(g)\}$ in the increasing order and lets the $i$-th round
last for $h(i)+1$ many moves, and \lib\ answers with sequences
$d(n)\in 2^{k_n}:n\in\bigcup\rng(g)$. The result $X$ of this play
contains the real $r$ by the choice of the functions $g$ and $d$. However, the play is in the ground model,
therefore the set $X$ is an $I$-small
closed set coded in the ground model and the generic real $r$ is forced to fall out of all such sets.
Contradiction!
\end{proof}

\section{Uniformity of the meager ideal}

The proof of Theorem~\ref{thirdtheorem} depends on a fusion game characterization of the property $\Theta(P_I)$ for
definable $\gs$-ideals $I$. Fix a partition $\gw=\bigcup_i a_i$
of the natural numbers into infinite sets and consider the game $G_\Theta(I)$ between players \old\ and \boz\ 
of length $\gw$. All moves in it are $I$-positive Borel sets again.
\boz\ starts out with an initial set $B_{ini}$ and then at each round $j$ \old\ plays a set $B_j$, which \boz\ answers
with its subset $C_j$. \old\ wins if the \emph{result} of the play, the set $B_{ini}\cap\bigcap_i\bigcup_{j\in a_i}C_j$,
does not belong to the ideal $I$.

\begin{lemma}
\label{thetaclaim}
(LC) $\Theta(P_I)$ if and only if \old\ has a winning strategy in the game $G_\Theta(I)$.
\end{lemma}

\begin{proof}
The right-to-left direction is easy. Let $\gs$ be a winning strategy for \old. The proof of properness
is the same as in Lemma~\ref{keylemma} and it is left to the reader. To see that the ground model
reals are not meager in the extension, let $B\in P_I$ be a condition and $\{\dot O_n:n\in\gw\}$ a name
for a sequence of open dense subsets of $2^{<\gw}$. Simulate a play of the game $G_\Theta(I)$
in which \boz\ plays the initial move $B=B_{ini}$ and then on the side constructs a chain
$t_0\subset t_1\subset\dots$ of finite binary sequences and plays so that for every integer $i\in\gw$
and every $j\in a_i$ the condition $C_j$ forces the sequence $\check t_j$ into $\dot O_n$. Let $r=\bigcup_nt_n$.
The result of the play is then a condition in the poset $P_I$ which forces the ground model real $\check r$
to have an initial segment in every set $\dot O_n:n\in\gw$. 

For the left-to-right direction 
note that the game is determined by the large cardinal assumptions, and it is enough to derive a contradiction
from the assumption that \boz\ has a winning strategy $\gs$. Towards the contradiction, choose a countable
elementary submodel $M$ of a large enough structure containing the strategy $\gs$ and let $T$ be the
tree of all partial plays of the game $G_\Theta(I)$ respecting the strategy $\gs$ in the model $M$ in which
\boz\ makes the last move. For a node $t\in T$ let $\last(t)$ be this last \boz's move in the play $t$.

It is clear that for every node $t\in T$ the set $D_t=\{\last(s):s\in T$ is an immediate successor of the node $t\}$
is the intersection of some dense set $E_t\subset P$ in the model $M$ with the model $M$ itself, namely of the dense set
$E_t=\{C\in P_I:\exists B\in P_I\ t^\smallfrown B^\smallfrown C$ respects the strategy $\gs\}$. Thus, if $B_{ini}\in M$
is the initial move dictated by the strategy $\gs$ and $p\leq B_{ini}$ is some $M$-master condition
below it, it is the case that $p\Vdash\forall t\in T\ \exists C\in D_t\ \dotrgen\in C$. Let $r$ be a generic real below the
condition $p$. Since the set of all ground model branches of the tree $T$ is not meager, 
there is a ground model branch $b\subset T$ such that for every number $i\in\gw$
there is $j\in a_i$ such that $r\in C_j$, where $C_j$ is the set the strategy $\gs$ played on the $j$-th round
of the rund $b$.

Now the play $b\subset T$ is in the ground model, and its resulting set is a ground model coded $I$-small Borel set.
The real $r$ belongs to the resulting set by the choice of the play $b$, but at the same time
it is forced to fall out of all such sets. Contradiction!
\end{proof}

There is an important corollary
which greatly simplifies certain statements in \cite{z:book}. Recall:

\begin{definition}
\cite{z:book}, Chapter 4.
A forcing $P$ is \emph{strongly proper} if for every countable elementary submodel $M$ of a large enough
structure, every condition $p_0\in P\cap M$ and every collection $\{D_i:i\in\gw\}$ of dense subsets of the poset
$P\cap M$ there is a \emph{strong master} condition $p_1\leq p_0$ 
forcing the generic filter to meet all the sets in the collection.
\end{definition}

\noindent It turns out that this notion is in the definable context identical to $\Theta$.

\begin{corollary}
(LC) Let $I$ be a definable $\gs$-ideal. $P_I$ is strongly proper if and only if $\Theta(P_I)$ holds.
\end{corollary}

\begin{proof}
First suppose that $\Theta(P_I)$ holds, and use Claim~\ref{thetaclaim} to find a winning strategy $\gs$ for \old\ in
the game $G_\Theta(P_I)$. Now suppose that $M$ is a countable elementary submodel of a large enough structure containing
$I$ and $\gs$, $B\in M\cap P_I$
is a condition, and $\{D_i:i\in\gw\}$ is a collection of dense subsets of the poset $P_I\cap M$.
Simulate a play of the game $G_\Theta(I)$ in which \old\ follows the strategy $\gs$ and \boz\ plays the
set $B$ as her initial move, and makes sure that for every number $i\in\gw$ and every $j\in a_i$, $C_j\in D_i$.
The resulting set of the play will be the desired strong master condition.

On the other hand, suppose that a poset $P$ is strongly proper. To show that the ground model reals are not meager,
just choose a condition $p_0\in P$ and a name $\{\dot O_i:i\in\gw\}$ for a collection of open dense
subsets of $2^{<\gw}$. Let $M$ be a countable elementary submodel of a large enough structure containing
all the relevant objects. Using the fact that the set $P\cap M$ is countable, it is easy to inductively construct
a sequence $t_0\subset t_1\subset\dots$ of finite binary sequences such that for all $i\in\gw$ and $p\in P\cap M$
there is some $j\in\gw$ and a condition $q\in M$, $q\leq p$ such that it forces $\check t_j\in\dot O_n$. 
This is to say that the sets $D_i=\{q\in P\cap M:\exists j\ q\Vdash t_j\in\dot O_i\}$ are dense in the
poset $P\cap M$ for all numbers $i\in\gw$. Let $p_1\leq p_0$ be the strong master condition, and let $r=
\bigcup_jt_j$. Clearly $p_1$ forces the real $\check r$ to have an initial segment in each of the sets $\dot O_n$ as desired.
\end{proof}

\noindent Towards the proof of Theorem~\ref{thirdtheorem}, let $I,J$ be definable $\gs$-ideals,
and suppose that $\Theta(I)$ and $\Phi(J)$ hold. 

\begin{corollary}
$\Theta(P_I\times P_J)$ holds.
\end{corollary}

\begin{proof}
It is enough to show that the poset is strongly proper.
Fix a winning strategy $\gs$ for \old\ in the game $G_\Phi(J)$, let $M$ be a countable
elementary submodel of a large enough structure, $p_0\times q_0\in M$ a Borel $I\times J$ block in it,
and $\{D_i:i\in\gw\}$ a collection of open dense subsets of the poset $P_I\times P_J\cap M$. 

First note that there is a $J$-positive Borel set $q_1\subset q_0$ such that for every number $i\in\gw$
the set $E_i=\{p\in P_I\cap M:\exists q\in M\ q_1\leq q\land \langle p,q\rangle\in D_i\}$ is open dense in $P_I\cap M$.
To see this, simulate a play of the game $G_\Phi(J)$ against the strategy $\gs$ in the same way as in the first
two paragraphs of the proof of Lemma~\ref{keylemma}, with the poset $2^{<\gw}$ replaced by $P_I\cap M$.

Now since the forcing $P_I$ is strongly proper, the set $p_1=p_0\cap\bigcap_i\bigcup E_i$ is Borel and $I$-positive.
It is not difficult to see that $p_1\times q_1\subset p_0\times q_0\cap\bigcap_i D_i$, and the block
$p_1\times q_1$ is the desired strong master condition. Thus the poset $P_I\times P_J$ is strongly proper.
\end{proof}

\begin{corollary}
$\RR(I,J)$ holds.
\end{corollary}

\begin{proof}
Suppose that $p_0\times q_0$ is an $I\times J$-block, decomposed into a countable union $\bigcup_mC_m$ of Borel sets.
It is enough to show that one of the sets $C_m:m\in\gw$ contains a block. Since $p_0\times q_0$ forces the generic
pair of reals into itself, there must be a strengthening $p_1\times q_1$ which forces the generic pair into
a set $C_m$ for some specific number $m\in\gw$. Let $M$ be a countable elementary submodel of a large enough structure
containing the sets $p_1, q_1, C_m$. The argument from the previous proof produces a block $p_2\times q_2\subset
p_1\times q_1$ consisting only of pairs of mutually $M$-generic reals. By the forcing theorem and an absoluteness
argument, $p_2\times q_2\subset C_m$ as desired.
\end{proof}

\noindent Theorem~\ref{thirdtheorem} follows.

\bibliographystyle{plain}
\bibliography{zapletal,odkazy,shelah}

\end{document}